# REGIONAL BOUNDARY GRADIENT STRATEGIC SENSORS IN NEUMANN BOUNDARY CONDITIONS

Raheam Al-Saphory[1*] and Mrooj Al-Bayati[1]

[1*] Department of Mathematics, College of Education For Pure Sciences, Tikrit University, Iraq.

[1*] Email: saphory@hotmail.com and mroojal.bayati@hotmail.com

*To whom correspondence should be addressed

**ABSTRACT**

The main idea of this paper is to explore and present the original results related to the notion of regional boundary gradient strategic sensors in distributed parameter system for Neumann problem. Thus, this systems is described by parabolic state space partial differential equations where the dynamic is governed by strongly continuous semi-group in Hilbert space. In addition, we give the sufficient conditions which are guaranteed the characterizations for such sensors in order that regional boundary approximately gradient observability can be achieved and proved. The obtained results are applied to two-dimensional circular in a part of considered system domain with various cases of sensors locations are tackled and studied.

## 1. INTRODUCTION

One of the most important concept in the analysis of distribution parameter systems is the observation problem [9-14]. In many situations one is interested in the knowledge of the state of a partial differential equations system on a sub-domain or region $\Gamma$ in the spatial boundary domain $\partial\Omega$ in finite time interval or infinite [1-8, 14-23].

The study of this notion motivated by certain concrete-real problem, in thermic, mechanic, environment [21-23]. This is the case, for example, of energy exchange problem, where the goal is to observe the energy exchanged between a casting plasma on a plane target which is perpendicular to the direction of the flow from measurements carried out by internal thermocouples (Figure 1).

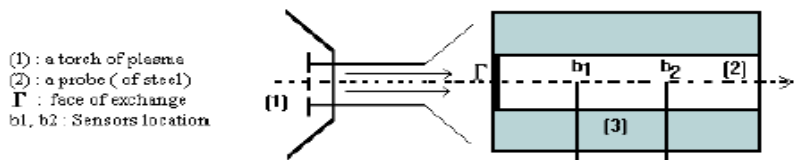

**Figure1.** The observation Profile of Energy Exchanged on $\Gamma$.

The aim of this study is to give sufficient conditions of the boundary gradient strategic sensors in this region which observes regional boundary gradient state. This paper is organized as follow.

Some definitions with characterizations of the regional boundary gradient observability are provided in the second section. The third section gives many of the sufficient conditions for the regional boundary gradient strategic sensors and a reconstruction method is developed to cross from internal region to the boundary. Finally, we present and illustrate some applications with many situations of sensor locations in circular domain.

## 2. $\Gamma_G$-OBSERVABILITY

This section is devoted to studies the concept of regional boundary gradient observability in Neumann boundary conditions and gives some important results which is related to this concept.

### 2.1 Considered System and Preliminaries

Let the following assumptions be given

- $\Omega$ be a regular bounded open subset of $R^n$, with smooth boundary $\partial\Omega$.

- $\Gamma$ be a sub-boundary of $\partial\Omega$.

- $[0,T]$, $T > 0$ be a time measurement interval.

- The Hilbert spaces $Z, U$ and $O$ are separable where $Z$ is the state space, $U = L^2(0,T,R^p)$ is the control space and $O = L^2(0,T,R^q)$ is the observation space, where $p$ and $q$ are the numbers of actuators and sensors [12].

- The considered parabolic distributed system described by the following equation

$$\begin{cases} \frac{\partial z}{\partial t}(\mu,t) = Az(\mu,t) + Bu(t) & \Omega \times ]0,T[ \\ z(\mu,0) = z_0(\mu) & \overline{\Omega} \\ \frac{\partial z}{\partial \nu}(\eta,t) = 0 & \partial\Omega \times ]0,T[ \end{cases} \quad (1)$$

where $Q = \Omega \times ]0,T[$, $\Sigma = \partial\Omega \times ]0,T[$, $\mu \in \Omega$, $\eta \in \partial\Omega$, $t \in [0,T]$ and $(\mu,t) \in \Omega \times ]0,T[$, $(\eta,t) \in \partial\Omega \times ]0,T[$, $(\mu,0) \in \overline{\Omega}$, where $\overline{\Omega}$ is the closure of $\Omega$, $\frac{\partial z}{\partial \nu}$ denote the normal derivative at the boundary. Thus, $A$ is a second order linear differential operator, which is generated a strongly continuous semi-group $(S_A(t))_{t \geq 0}$ on the Hilbert space $Z = H^1(\Omega)$ and is self-adjoint with compact resolvent. The system (1) is augmented with the output function

$$y(.,t) = Cz(.,t) \quad (2)$$

and then:

- The operators $B \in L(R^p, Z)$ and $C \in L(Z, R^q)$ depend on the structure of actuators and sensors [2]. Consequently, the system (1) under the given assumptions has a unique solution [9-10] given as follows

$$z(\mu,t) = S_A(t)z_0(\mu) + \int_0^t S_A(t-\tau)Bu(\tau)d\tau \quad (3)$$

- The problem is how to give sufficient conditions of regional boundary gradient strategic sensors in a given sub-boundary $\Gamma$ in Neumann boundary conditions.

- Now, we define the operator

$$K: z \in Z \to Kz = CS_A(.)z \in O$$

and, the adjoint operator of $K$ denotes by $K^*$ define by

$$K^*y^* = \int_0^t S_A^*(s)C^*y^*(s)ds$$

- Consider the gradient operator

$$\begin{cases} \nabla: H^1(\Omega) \to (H^1(\Omega))^n \\ z \to \nabla z = (\frac{\partial z}{\partial \mu_1}, \ldots, \frac{\partial z}{\partial \mu_n}) \end{cases}$$

and the adjoint of $\nabla$ denotes by $\nabla^*$ is given as

$$\begin{cases} \nabla^*: (H^1(\Omega))^n \to H^1(\Omega) \\ z \to \nabla_z^* = u \end{cases}$$

where $u$ is a solution of the Neumann problem

$$\begin{cases} \Delta u = -f(z) & \text{in } \Omega \\ \partial u/\partial \nu = 0 & \text{in } \partial\Omega \end{cases}$$

- The trace operator of order zero is given by

$$\gamma_0: H^1(\Omega) \to H^{1/2}(\partial\Omega)$$

Thus, the extension of the trace operator which is denoted by $\gamma$ is linear, subjective and continuous, defined as

$$\gamma : (H^1(\Omega))^n \to (H^{1/2}(\partial\Omega))^n$$

and the adjoints are respectively given by $\gamma_0^*, \gamma^*$.

- For a sub-boundary $\Gamma \subset \partial\Omega$, we consider a gradient restriction operator

$$\chi_\Gamma : (H^{1/2}(\partial\Omega))^n \to (H^{1/2}(\Gamma))^n$$

and

$$\tilde{\chi}_\Gamma : H^{1/2}(\partial\Omega) \to H^{1/2}(\Gamma)$$

where the adjoints are respectively given by $\chi_\Gamma^*, \tilde{\chi}_\Gamma^*$.

- For a non-empty subset $\omega$ of $\Omega$, let $\chi_\omega$ be the restriction function defined by

$$\chi_\omega : \begin{cases} (H^1(\Omega))^n \to (H^1(\omega))^n \\ z \to \chi_\omega z = z|_\omega \end{cases}$$

where $z|_\omega$ is the restriction of the state $z$ to $\omega$ [7]. There adjoints are respectively denoted by $\chi_\omega^*$ are defined by

$$\chi_\omega^* : \begin{cases} (H^1(\omega))^n \to (H^1(\Omega))^n \\ z \to \chi_\omega^* z = \begin{cases} z|_\omega & \text{in } \omega \\ 0 & \text{in } \Omega \setminus \omega \end{cases} \end{cases}$$

- Finally, we denote the operators $H_\Gamma = \chi_\Gamma \gamma \nabla K^*$ from $\mathcal{O}$ into $(H^{1/2}(\Gamma))^n$ and $H_\omega = \chi_\omega \nabla K^*$ from $\mathcal{O}$ into $(H^1(\omega))^n$.

**2.2 Definitions and Characterizations**

In this sub-section, we present some definitions and results about the regional boundary gradient observability and strategic sensors. For that purpose, consider the autonomous system of (1) is defined by the following parabolic equatios

$$\begin{cases} \frac{\partial z}{\partial t}(\mu, t) = Az(\mu, t) & \Omega \times ]0, T[ \\ z(\mu, 0) = z_0(\mu) & \bar{\Omega} \\ \frac{\partial z}{\partial v}(\eta, t) = 0 & \partial\Omega \times ]0, T[ \end{cases} \qquad (4)$$

The solution of (4) is given by the following form
$$z(\mu, t) = S_A(t) z_0(\mu) \quad \text{for all } t \in [0, T] \qquad (5)$$

**Definition 2.1:** The system (4) augmented with the output function (2) (or the system (1)-(2)) is said to be an exactly (respectively approximately) regionally gradient observable on $\omega$, if
$$\text{Im } H_\omega = (H^1(\omega))^n \quad (\text{respectively } \overline{\text{Im } H_\omega} = (H^1(\omega))^n)$$

**Definition 2.2:** The system (4)-(2) is said to be an exactly (approximately respectively) regionally boundary gradient observable on $\Gamma$, if
$$\text{Im } H_\Gamma = (H^{1/2}(\Gamma))^n \quad (\text{respectively } \overline{\text{Im } H_\Gamma} = (H^{1/2}(\Gamma))^n)$$

**Remark 2.3:** We deduced that, this equation $\overline{\text{Im } H_\Gamma} = (H^{1/2}(\Gamma))^n$ is equivalent to $\ker H_\Gamma^* = \{0\}$.

**Proposition 2.4:** The system (4)-(2) is exactly $\Gamma_G$-observable if and only if there exists $v > 0$, such that for all $z^* \in (H^{1/2}(\Gamma))^n$,

$$\|\chi_\Gamma z^*\|_{(H^{1/2}(\Gamma))^n} \le v \|K\nabla^* \gamma^* \chi_\Gamma^* z^*\|_{\mathcal{O}} \qquad (6)$$

**Proof:** The proof of this proposition is conclude from the following general result as in [9]. Let $E, F$ and $G$ be a reflexive Banach spaces and $f \in L(E, G)$, $g \in L(F, G)$, then the following properties are equivalent

1. $Imf \subset Img$.

2. There exists $v > 0$, such that

$$\|f^* z^*\|_{E^*} \leq v \|g^* z^*\|_{F^*}, \text{ for all } z^* \in G^*.$$

If we apply this result, considered $E = G = (H^{1/2}(\Gamma))^n, F = O, f = Id_{(H^{1/2}(\Gamma))^n}$ and $g = \chi_\Gamma \gamma \nabla K^*$. Therefore, we obtain the inequality

$$\|\chi_\Gamma z^*\|_{(H^{1/2}(\Gamma))^n} \leq v \|K\nabla^* \gamma^* \chi_\Gamma^* z^*\|_o. \square$$

Now, we can get the following proposition.

**Proposition 2.5:** The system is exactly regionally boundary gradient observable (or exactly $\Gamma_G$-observable), if it is exactly boundary gradient observable (exactly $B_G$-observable).

**Proof:** Since the system is an exactly $B_G$-observable. Then there exists $\gamma > 0$, such that for all $z_0 \in H^{1/2}(\partial\Omega)$, we have

$$\|z_0\|_{H^{1/2}(\Gamma)} \leq \gamma \|K\gamma_0^* \tilde{\chi}_\Gamma^* z_0\|_{L^2(0,T,o)}, \text{for all } \gamma > 0$$

since $(H^{1/2}(\Gamma))^n \subset H^{1/2}(\Gamma)$, then

$$\|\nabla z_0\|_{(H^{1/2}(\partial\Omega))^n} = \|z_0\|_{(H^{1/2}(\Gamma))^n} \leq \|z_0\|_{H^{1/2}(\Gamma)}, \text{ for all } z_0 \in H^{1/2}(\Gamma)$$

where,

$$H^{1/2}(\Gamma) = \{z_0 : \int_\Gamma |z_0|^2 < \infty\}$$

and,

$$(H^{1/2}(\Gamma))^n = \left\{\nabla z_0 = g_i : \int_\Gamma |g_i|^2 < \infty, g_i = \frac{\partial z_0}{\partial \mu_i}, \text{for all } i = 1, 2, \dots\right\} \quad (7)$$

to prove that $\|z_0\|_{(H^{1/2}(\Gamma))^n} \leq v \|K\nabla^* \gamma^* \chi_\Gamma^* z_0\|_{L^2(0,T,o)}$

Now from using (7) and since the system is exactly $B_G$-observable, then there exists $\gamma > 0$ and $v > 0$, This allow to choose that $\gamma = \frac{1}{v}$ and by setting

$$v = \frac{\|K \gamma_0^* \tilde{\chi}_\Gamma^* z_0\|_o}{\|K \nabla^* \gamma^* \chi_\Gamma^* z_0\|_o} \quad (8)$$

and then,

$$\|z_0\|_{(H^{1/2}(\Gamma))^n} \leq \|z_0\|_{H^{1/2}(\Gamma)} \leq \gamma \|K \gamma_0^* \tilde{\chi}_\Gamma^* z_0\|_o \quad (9)$$

substitute (8) in (9), we obtain

$$\|z_0\|_{(H^{1/2}(\Gamma))^n} \leq \|K\nabla^* \gamma^* \chi_\Gamma^* z_0\|_o.$$

Therefore, this system is an exactly $\Gamma_G$-observable with $\gamma = 1$. $\square$

## 3. SUFICIENT CONDITIONS FOR $\Gamma_G$-STRATEGIC SENSORS

For achieving the regional boundary gradient observability, we will give the sufficient condition for the characterizations of sensors in a given sub-region $\Gamma$.

### 3.1 Concept of Sensors

We recall that some definitions and types of sensors as in [11].

**Definition 3.1:** A sensor may be defined by any couple $(D, f)$, where $D$ is spatial support represented by a non-empty part of $\overline{\Omega}$ and $f$ represent the distribution of the sensing measurements on $D$.

According to the choice of the parameters $D$ and $f$, we could have different types of sensors.

• A sensor may be point wise when $D = \{b\}$, with $b \in \overline{\Omega}$ and $f = \delta(.-b)$, where $\delta$ is the Dirac mass concentrated in $b$. In this case, the operator $C$ is unbounded and the output function (2) may be given by the form as in [9-10]

$$y(t) = \int_\Omega z(\mu, t) \delta_b(\mu - b) d\mu = z(b, t)$$

• A sensor may be zone, if $D \subset \Omega$ and $f \in L^2(D)$. The output function (2) may be given by the form

$$y(t) = \int_D z(\mu, t) f(\mu) d\mu.$$

• It may be boundary zone if $\Gamma_i \subset \partial\Omega$ and $f_i \in L^2(\Gamma_i)$, the output function (2) may be given by the form

$$y(t) = \int_{\Gamma_i} z(\eta, t) f_i(\eta) d\eta$$

**Definition 3.2:** A sensor $(D, f)$ is $\Gamma_G$-strategic, if the observed system are approximately $\Gamma_G$-observable.

**Definition 3.3:** A suit of sensors $(D_i, f_i)_{1 \le i \le q}$ are $\Gamma_G$-strategic, if there exist at least one sensor $(D_1, f_1)$ which is $\Gamma_G$-strategic.

**Proposition 3.4:** A sensor is $\Gamma_G$-strategic if and only if the operator $N_\Gamma = HH^*$ is positive definite.

**Proof:** Since a sensor is $\Gamma_G$-strategic this mean that the system is approximately $\Gamma_G$-observable, let $z^* \in (H^{1/2}(\Gamma))^n$, such that

$\langle N_\Gamma z^*, z^* \rangle_{(H^{1/2}(\Gamma))^n} = 0$, then $\|H^* z^*\|_0 = 0$

and therefore $z^* \in \ker H^*$, thus $z^* = 0$, i.e., $N_\Gamma$ is positive definite.

Conversely, let $z^* \in (H^{1/2}(\Gamma))^n$, such that

$H^* z^* = 0$, then $\langle H^* z^*, H^* z^* \rangle_0 = 0$

and thus,

$\langle N_\Gamma z^*, z^* \rangle_{(H^{1/2}(\Gamma))^n} = 0.$

Hence $z^* = 0$, thus the system is approximately $\Gamma_G$-observable and therefore a sensor is $\Gamma_G$-strategic. □

**Proposition 3.5:** A sensor is $\Gamma_G$-strategic, if the observed systems are an exactly $\Gamma_G$-observable.

**Proof:** Since the system is an exactly $\Gamma_G$-observable. Then, we have

$\operatorname{Im} H_\Gamma = (H^{1/2}(\Gamma))^n$

from decomposition sub-spaces of direct sum in a Hilbert space, we represent $(H^{1/2}(\partial\Omega))^n$ by the unique form

$$\ker \chi_\Gamma + Im\, \chi_\Gamma^* \chi_\Gamma \gamma \nabla K^* = (H^{1/2}(\partial\Omega))^n$$

we obtain that,    $\ker K(t) \nabla^* \gamma^* \chi_\Gamma^* = \{0\}$

and this is equivalent to    $\overline{Im\, \chi_\Gamma \gamma \nabla K^*} = (H^{1/2}(\Gamma))^n$.

Finally, we can deduced this system is approximately $\Gamma_G$-observable and then this sensor is $\Gamma_G$-strategic. □

**Remark 3.6:** From the previous results, we can deduced that.

1. If the systems are an exactly $\Gamma_G$-observable, then the systems are approximately $\Gamma_G$-observable, and therefore the sensors are $\Gamma_G$-strategic.

2. A sensor which is a regional boundary gradient strategic in $\Gamma_G^1$ ($\Gamma_G^1$-strategic sensor) for a system, where $\Gamma_G^1 \subset \Gamma_G$ is a regional boundary gradient strategic sensor in $\Gamma_G^2$ ($\Gamma_G^2$-strategic sensor) for any $\Gamma_G^2 \subset \Gamma_G^1$.

**3.2 The Main Result**

This sub-section is related to develop the results which are concerned with the regional boundary gradient strategic sensors, and give the sufficient conditions for each sensor. For this purpose we assume that there exists a complete set of eigenfunctions $(\varphi_{nj})_{n \in I, j=1,\ldots,m_n}$ of $A$ in $H^1(\bar{\Omega})$ associated with eigenvalue $\lambda_n$ of multiplicities $m_n$ and $m_n = \sup_{n \in I} m_n$ is finite. For $\bar{z} = (z_1, \ldots, z_{n-1})$ and $\bar{n} = (n_1, \ldots, n_{n-1})$. Suppose that the function $\psi_{\bar{n}j}(\bar{z}) = \chi_\Gamma \gamma \nabla \varphi_{nj}(z)$, $n \in I$, is a complete set in $(H^{1/2}(\Gamma))^n$. If the system (1) has $J$ unstable modes, then we have the following theorem.

**Theorem 3.7:** suppose that $\sup m_n = m < \infty$, then the suite of sensors $(D_i, f_i)_{1 \le i \le q}$ are $\Gamma_G$-strategic sensor if and only if

1. $q \ge m$,
2. $\operatorname{rank} G_n = m_n$, for all $n \ge 1$, where $G_n = (G_n)_{ij}$ with $1 \le i \le q, 1 \le j \le m_n$, and

$$(G_n)_{ij} = \begin{cases} \sum_{k=1}^n \frac{\partial \psi_{\bar{n}j}}{\partial z_k}(b_i) & \text{point wise sensor} \\ \sum_{k=1}^n \langle \frac{\partial \psi_{\bar{n}j}}{\partial z_k}, f_i \rangle_{L^2(D_i)} & \text{zone sensor} \end{cases}$$

**Proof:** First, we recall that the systems are approximately $\Gamma_G$-observable are equivalent to $[K\nabla^* \gamma^* \chi_\Gamma^* z^* = 0 \Rightarrow z^* = 0]$ which allows to say that the sequence of sensors $(D_i, f_i)_{1 \le i \le q}$ are $\Gamma_G$-strategic if and only if

$$\{z^* \in (H^{1/2}(\Gamma))^n | \langle Hy, z^* \rangle_{(H^{1/2}(\Gamma))^n} = 0, \text{for all } y \in O\} = \{0\}$$

suppose that the suite of sensors $(D_i, f_i)_{1 \le i \le q}$ are $\Gamma_G$-strategic on $\Gamma$, but for a certain $n \in N$, $\operatorname{rank} G_n \ne m_n$, then there exists a vector $z_n = (z_{n_1}, z_{n_2}, \ldots, z_{n_m})^T \ne 0$, such that $G_n z_n = 0$.

so, we can construct a non-zero $z_0 \in H^{1/2}(\Gamma)$ considering $\langle z_0, \psi_{pj} \rangle_{H^{1/2}(\Gamma)} = 0$, if $p \ne n$, and $\langle z_0, \psi_{nj} \rangle_{H^{1/2}(\Gamma)} = z_{nj}$, $1 \le j \le m_n$, $z_0 = \sum_{j=1}^{m_n} z_{nj} \psi_{nj} \in H^{1/2}(\Gamma) \ne 0$, then

$$\langle Hy, z_0 \rangle_{(H^{1/2}(\Gamma))^n} = \sum_{k=1}^n \langle \tilde{\chi}_\Gamma \gamma_0 \frac{\partial}{\partial \mu_k}(K^* y), \tilde{\chi}_\Gamma^* z_0 \rangle_{H^{1/2}(\Gamma)}$$

$$= \sum_{k=1}^n \langle \frac{\partial}{\partial \mu_k}(\tilde{z}(T)), \gamma_0 \tilde{\chi}_\Gamma^* z_0 \rangle_{H^{1/2}(\partial\Omega)}$$

where $\tilde{z}$ is the solution of the following system

$$\begin{cases} \frac{\partial \tilde{z}}{\partial t}(\mu, t) = A^*\tilde{z}(\mu,t) + \sum_{i=1}^{q} f_i\, y_i(T-t) & \Omega \times ]0,T[ \\ \tilde{z}(\mu,0) = 0 & \bar{\Omega} \\ \frac{\partial \tilde{z}}{\partial v}(\eta, t) = 0 & \partial\Omega \times ]0,T[ \end{cases} \quad (10)$$

Now, consider the system

$$\begin{cases} \frac{\partial \psi}{\partial t}(\mu, t) = -A\varphi(\mu,t) & \Omega \times ]0,T[ \\ \psi(\mu,0) = \gamma_0^* \tilde{\chi}_T^* z_0 & \bar{\Omega} \\ \frac{\partial \psi}{\partial v}(\eta, t) = 0 & \partial\Omega \times ]0,T[ \end{cases} \quad (11)$$

multiply (10) by $\frac{\partial \psi}{\partial \mu_k}$ and integrate on $Q$, we obtain that

$$\int_Q \frac{\partial \psi}{\partial \mu_k}(\mu, t) \frac{\partial \tilde{z}}{\partial t}(\mu, t)\, d\mu\, dt = \int_Q A^*\tilde{z}(\mu, t) \frac{\partial \psi}{\partial \mu_k}(\mu, t)\, d\mu\, dt$$
$$+ \int_Q \left(\sum_{i=1}^{q} \delta_{b_i} y_i(T-t)\right) \frac{\partial \psi}{\partial \mu_k}(\mu, t)\, d\mu\, dt$$

But, we have

$$\int_Q \frac{\partial \varphi}{\partial \mu_k}(\mu, t) \frac{\partial \tilde{z}}{\partial t}(\mu, t)\, d\mu\, dt = \int_{\partial\Omega} \left[\frac{\partial \varphi}{\partial \mu_k}(\mu, t)\tilde{z}(\mu, t)\, d\mu\right]_0^T +$$

$$\int_Q A\frac{\partial \psi}{\partial \mu_k}(\mu, t)\tilde{z}(\mu, t)\, d\mu\, dt = \int_{\partial\Omega} \frac{\partial \psi}{\partial \mu_k}(\mu, t)\tilde{z}(\mu, t)\, d\mu + \int_Q A\frac{\partial \psi}{\partial \mu_k}(\mu, t)\tilde{z}(\mu, t)\, d\mu\, dt.$$

then,

$$\int_{\partial\Omega} \frac{\partial \psi}{\partial \mu_k}(\mu, t)\tilde{z}(\mu, t)\, d\mu = -\int_Q A\frac{\partial \psi}{\partial \mu_k}(\mu, t)\tilde{z}(\mu, t)\, d\mu$$
$$+ \int_Q A^*\tilde{z}(\mu, t)\frac{\partial \psi}{\partial \mu_k}(\mu, t)\, d\mu\, dt + \int_Q \left(\sum_{i=1}^{q} \delta_{b_i} y_i(T-t)\right) \frac{\partial \psi}{\partial \mu_k}(\mu, t)\, d\mu\, dt.$$

Integrating by parts, we obtain

$$\int_{\partial\Omega} \frac{\partial \psi}{\partial \mu_k}(\mu, t)\tilde{z}(\mu, t)\, d\mu = -\int_\pi \frac{\partial \tilde{z}(\eta, t)}{\partial v_{A^*}} \frac{\partial \psi}{\partial \mu_k}(\eta, t)\, d + \int_\pi \frac{\partial}{\partial v_{A^*}}\left(\frac{\partial \psi}{\partial \mu_k}(\eta, t)\, d\eta dt\right)\tilde{z}(\eta, t)$$
$$+ \int_Q \left(\sum_{i=1}^{q} \delta_{b_i} y_i(T-t)\right) \frac{\partial \psi}{\partial \mu_k}(\mu, t)\, d\mu\, dt.$$

the boundary conditions give

$$\int_{\partial\Omega} \frac{\partial \psi}{\partial \mu_k}(\mu, t)\tilde{z}(\mu, t)\, d\mu = \int_Q \left(\sum_{i=1}^{q} \delta_{b_i} y_i(T-t)\right) \frac{\partial \psi}{\partial \mu_k}(\mu, t)\, d\mu\, dt.$$

Thus,

$$\int_{\partial\Omega} \psi(\mu, t) \frac{\partial \tilde{z}}{\partial \mu_k}(\mu, T)\, d\mu = -\sum_{i=1}^{q} \int_0^T \frac{\partial \psi}{\partial \mu_k}(b_i, t)\, y_i(T-t)\, dt.$$

and, we have

$$\langle \chi_\Gamma \gamma \nabla K^*, z_0 \rangle_{(H^{1/2}(\Gamma))^n} = \sum_{k=1}^n \int_\Omega \frac{\partial \tilde{z}}{\partial \mu_k}(\mu, t) \psi(\mu, t) d\mu$$

$$= -\sum_{k=1}^q \int_0^T \sum_{k=1}^n \frac{\partial \psi}{\partial \mu_k}(b_i, t) y_i(T-t) dt$$

but,

$$\psi(\mu, t) = \sum_{p=1}^\infty e^{-\lambda_p(T-t)} \sum_{j=1}^{m_p} \langle z_0, \psi_{pj} \rangle_{L^2(\omega)} \psi_{pj}$$

Then,

$$\sum_{k=1}^n \frac{\partial \psi}{\partial \mu_k}(b_i, t) = \sum_{p=1}^\infty e^{-\lambda_p(T-t)} \sum_{j=1}^{m_p} \langle z_0, \varphi_{pj} \rangle_{L^2(\omega)} \sum_{k=1}^n \frac{\partial \psi}{\partial \mu_k}(b_i)$$

$$= \sum_{p=1}^\infty e^{\lambda_p(T-t)}(G_p z_p)_i$$

therefore,

$$\langle \chi_\Gamma \gamma \nabla K^* y, z_0 \rangle_{(H^{1/2}(\Gamma))^n} = -\sum_{K=1}^q \int_0^T \sum_{p=1}^\infty e^{\lambda_p(T-t)}(G_p z_p) y_i(T-t) dt$$

thus,

$$\langle \chi_\Gamma \gamma \nabla K^* y, z_0 \rangle_{(H^{1/2}(\Gamma))^n} = -\sum_{i=1}^q \int_0^T e^{\lambda_n(T-t)}(G_n z_n)_i \, y_i(T-t) dt = 0.$$

This is true for all $y \in L^2(0,T,R^q)$, then $z_0 \in \ker H_\Gamma^*$ which contradicts the assumption that the suite of sensors are $\Gamma_G$-strategic.□

**3.3 Internal and Boundary Reconstruction State Gradient Method**

The regional boundary gradient strategic sensors in $\Gamma$ ($\Gamma_G$-strategic sensor) may be seen as internal regional gradient strategic sensors if we consider $\bar{\omega}_r \subset \bar{\Omega}$ [23].

• The operator $\mathcal{R}$ is a continuous and linear extension which is given by

$$\mathcal{R}: (H^{1/2}(\partial\Omega))^n \to (H^1(\Omega))^n, \text{ such that}$$

$$\gamma \nabla \mathcal{R} h(\mu, t) = h(\mu, t), \quad \text{for all } h(\mu, t) \in (H^{1/2}(\partial\Omega))^n \qquad (12)$$

• Let $r > 0$ be an arbitrary and sufficiently small real number and let $E = \bigcup_{z \in \Gamma} B(z, r)$ and $\bar{\omega}_r = E \cap \bar{\Omega}$, where $B(z, r)$ is a ball of radius $r$ centered in $z(\mu, r)$, and $\Gamma$ is a part of $\bar{\omega}_r$ (Figure 2).

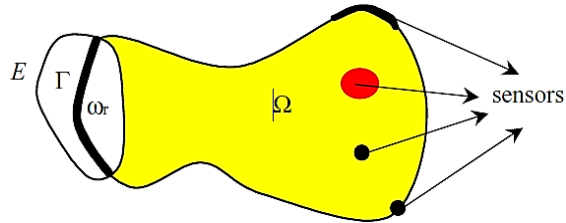

**Figure 2:** The Domain $\Omega$, region $\bar{\omega}_r$ and the boundary region $\Gamma$.

Now, we show that relation between the $\Gamma_G$-strategic sensor and $\bar{\omega}_{r_G}$-strategic sensor in the following result.

**Proposition 4.2:** From the above results, we deduce that

1. A sensor is $\Gamma_G$-strategic if the sensor is $\bar{\omega}_{r_G}$-strategic.

2. A sensor is $\Gamma_G$-strategic if the system is exactly $\bar{\omega}_{r_G}$-observable.

**Proof: 1.** Since the sensor is $\bar{\omega}_{r_G}$-strategic in $\bar{\omega}_r$, this mean that the system is approximately $\bar{\omega}_{r_G}$-observable in $\bar{\omega}_r$.

Thus, the system is approximately $\Gamma_G$-observable (Zerrik and Bourray 2003). Therefore, the sensor is $\Gamma_G$-strategic.□

**2.** Let $z(\mu, t) \in (H^{1/2}(\Gamma))^n$ and $\bar{z}(\mu, t)$ be an extension to $(H^{1/2}(\partial\Omega))^n$. By using equation (12) and trace theorem there exists $\mathcal{R}\bar{z}(\mu, t) \in (H^1(\Omega))^n$, with bounded support such that

$$\gamma \mathcal{R}\bar{z}(\mu, t) = \bar{z}(\mu, t).$$

Since the system is an exactly $\bar{\omega}_{r_G}$-observable, then the system is approximately $\bar{\omega}_{r_G}$-observable [3]. And, since a system is approximately $\bar{\omega}_{r_G}$-observable then a system is approximately $\Gamma_G$-observable [1, 21-23] Thus, the sensor is $\Gamma_G$-strategic.□

**Remark 4.3:** From the previous results, we have.

1. If the system is an exactly $\bar{\omega}_{r_G}$-observable, then the system is an exactly $\Gamma_G$-observable, i.e., there exists an operator $\chi_{\bar{\omega}_r} \nabla K^* : \mathcal{O} \to (H^1(\omega_r))^n$ given by

$$H_{\bar{\omega}_r} y(.,t) = \chi_{\bar{\omega}_r} \nabla K^* y(.,t) = \chi_{\bar{\omega}_r} \mathcal{R}\bar{z}(\mu, t).$$

Hence, $\chi_\Gamma \left( \gamma \chi_{\bar{\omega}_r} \nabla K^* y(.,t) \right) = z(\mu, t).$

where $z(\mu, t) \in (H^{1/2}(\Gamma))^n$ and $\bar{z}(\mu, t)$ be an extension to $(H^{1/2}(\partial\Omega))^n$.

2. If the system is approximately $\bar{\omega}_{r_G}$-observable, then the system is approximately $\Gamma_G$-observable.

3. An extension of these results can be applied for different cases of regional exponential reduced observability [4] and, to the regional boundary exponential reduced observability in distributed parameter systems [6].

## 4. APPLICATIONS TO SENSORS LOCATIONS IN CIRCULAR DOMAIN

In this section we study and explore the following cases:

**4.1 Zone Sensor Case**

In this sub-section, we shall discuss the characterization of sensors in the case of (internal and boundary) zone sensors in in region of circular boundary domain with Neumann boundary conditions.

**4.1.1 Internal circular zone sensor**

In this case, we consider the system given by the following form

$$\begin{cases} \frac{\partial z}{\partial t}(r, \theta, t) = \frac{\partial^2 z}{\partial r^2}(r, \theta, t) + \frac{\partial^2 z}{\partial \theta^2}(r, \theta, t) + z(r, \theta, t) & \Omega \times ]0, T[ \\ z(r, \theta, 0) = z_0(r, \theta) & \bar{\Omega} \\ \frac{\partial z}{\partial v}(a, \theta, t) = 0 & \partial\Omega \times ]0, T[ \end{cases} \quad (13)$$

where, $\Omega = ]0, a[$, $r = a > 0$, $\theta \in [0, 2\pi]$ are defined as in (Figure 3). The augmented output function is defined by

$$y_n(t) = \int_{D_n} \frac{\partial z}{\partial v}(r_n, \theta_n, t) f(r_n, \theta_n) dr_n d\theta_n, \quad 2 \leq n \leq q \tag{14}$$

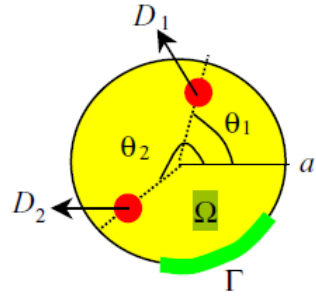

**Figure 3:** Circular domain $\Omega$, region $\Gamma$ and location $D_1, D_2$ of internal zone sensors.

So, the eigenfunctions and eigenvalues concerning the region $\Gamma = (a, \theta_i)_{2 \leq i \leq q}$ of $\partial \Omega$ with $\theta \in [0, 2\pi]$ be defined

$$\lambda_{nm} = -\beta_{nm}^2, \quad n \geq 0, m \geq 1 \tag{15}$$

where $\beta_{nm}$ are the zones of the Bessel functions $J_n$ and

$$\begin{cases} \varphi_{0m}(r, \theta) = m_0(\beta_{0m} r), \quad m \geq 1 \\ \varphi_{nm_1}(r, \theta) = m_n(\beta_{nm_1} r)\cos(n\theta), \quad n, m_1 \geq 1 \\ \varphi_{nm_2}(r, \theta) = m_n(\beta_{nm_2} r)\sin(n\theta), \quad n, m_2 \geq 1 \end{cases} \tag{16}$$

with multiplicity $m_n = 2$, for all $n, m \neq 0$ and $m_n = 1$, for all $n, m = 0$ as in [12]. In this case, there exist at least two zone sensors $(D_n, f_n)_{2 \leq n \leq q}$ where $D_n = (r_n, \theta_n)_{2 \leq n \leq q}$. Thus, we have the result.

**Proposition 4.1:** Suppose that $f_n$ and $D_n$ are symmetric with respect to $\theta = \theta_n$, for all $n, 2 \leq n \leq q$. Then the sensors $(D_n, f_n)$ are $\Gamma_G$-strategic, if

$$\frac{n_0(\theta_1 - \theta_2)}{\pi} \notin Q, \text{ for some } 1 \leq n_0 \leq J.$$

**4.1.2 Boundary zone sensor**

We consider the system (13) with the output function describe by

$$y(t) = \int_{\Gamma_n} z(r, \theta_n, t) f(r, \theta_n) d\theta_n, \quad \theta_n \in [0, 2\pi], \ t > 0 \tag{17}$$

In this case, it is necessary to have at least two boundary zone sensors $(\Gamma_n, f_n)_{2 \leq n \leq q}$ where $\Gamma_n = (r, \theta_n)_{2 \leq n \leq q}$ as in (Figure 4).

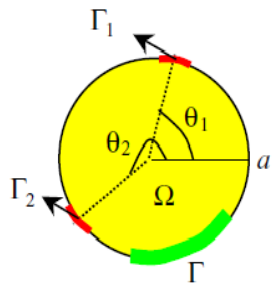

**Figure 4:** Circular domain $\Omega$, region $\Gamma$ and location $\Gamma_1, \Gamma_2$ of boundary zone sensors.

Thus, we have the result

**Proposition 4.2:** Suppose that $f_n$, $\Gamma_n$ are symmetric with respect to $\theta = (\theta_n)_{2\leq n\leq q}$. Then the sensors $(\Gamma_n, f_n)_{2\leq n\leq q}$ are $\Gamma_G$-strategic, if

$$\frac{n_0(\theta_1-\theta_2)}{\pi} \notin Q, \text{ for some } 1 \leq n_0 \leq J.$$

### 4.2 Case of Pointwise Sensor

In this sub-section, we discuss and describe the sensors in the case of (internal and boundary) point wise sensors in circular domain.

#### 4.2.1 Internal pointwise sensor

Consider the system (13) with the output function defined by the form

$$y_n(t) = \int_\Omega \frac{\partial z}{\partial v}(r_n, \theta_n, t) f(r_n, \theta_n) dr_n d\theta_n \qquad (18)$$

where $2 \leq n \leq q$, $\theta_n \in [0, 2\pi]$, $r = a > 0$, and $\Omega = (0, a)$, $t > 0$. The sensors may be located in $p_1 = (r_1, \theta_1)$ and $p_2 = (r_2, \theta_2) \in \Omega$ as in (Figure 5).

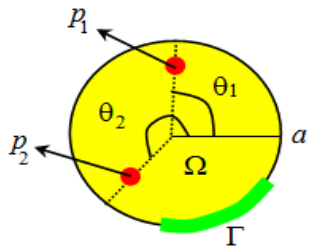

**Figure 5:** Circular domain $\Omega$, region $\Gamma$ and locations $p_1, p_2$ of internal pointwise sensors.

**Proposition 4.3:** We deduced that:

1. The sensors $(b, \delta_b)$ are $\Gamma_G$-strategic, if $\frac{n_0(\theta_1-\theta_2)}{\pi} \notin Q$, for all $1 \leq n_0 \leq J$.
2. If $r_1 = r_2$, then also, the sensors $(b, \delta_b)$ are $\Gamma_G$-strategic, if

$$\frac{n_0(\theta_1-\theta_2)}{\pi} \notin Q, \text{ for all } 1 \leq n_0 \leq J.$$

#### 4.2.2 Boundary pointwise sensor

In this case, we consider the system (16) with the output function is given by

$$y(t) = \int_{\partial\Omega} z(a, \theta_n, t) f(a, \theta_n) d\theta_n, \quad \theta_n \in [0, 2\pi], \quad 2 \leq n \leq q, \quad t > 0 \qquad (19)$$

When the pointwise sensors at the polar coordinates $p_n = (1, \theta_n)$ where $\theta_n \in [0, 2\pi]$, $2 \leq n \leq q$ and $t > 0$, as in (Figure 6).

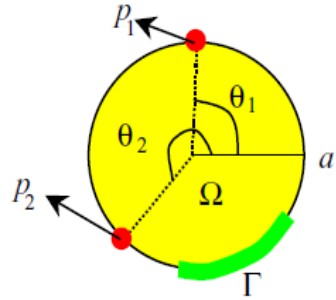

**Figure 6:** Circular domain, regional $\Gamma$ and locations $p_1, p_2$ of boundary point wise sensors.

Thus, we have the following result.

**Proposition 4.4:** The sensors $(b, \delta_b)$ are $\Gamma_G$-strategic, if $\frac{n_0(\theta_1 - \theta_2)}{\pi} \notin Q$, for all $1 \leq n_0 \leq J$.

**Remark 4.5:** These results can be extended to the following
1. Case of Dirichlet boundary conditions [5].
2. Case of rectangular domain in various case of point wise zone internal or boundary as in [3,5].
3. The results show that there are some sensor locations to be avoided [2].

**Acknowledgements.** Our thanks in advance to the editors and experts for considering this paper to publish in this esteemed journal. The authors appreciate your time and effort in reviewing the manuscript and greatly value your assistant as reviewer for the paper.

**5. CONCLUSION**

We have shown that the regional boundary gradient strategic sensors notion in Neumann problem, is guaranteed the existing of approximately regional boundary gradient observability concept in a region $\Gamma$. For infinite dimensional systems in Hilbert space, many interesting results concerning the choice of sensors structure are given and illustrated in specific situations. Finally, we have been shown that there is a relation between the regional boundary gradient strategic sensor with number and shapes sensors and the system domain. Many problems is stilled opened for example, the possibility to extend these results for the case of regional gradient detectability in hyperbolic systems [8].